\documentclass[11pt, english]{article}
\usepackage[margin=3cm]{geometry}
\usepackage{amsbsy,amsmath,amsthm,amssymb}

\usepackage{fourier}

\usepackage{xcolor,bm}
\usepackage{esint}
\usepackage{babel}
\usepackage[colorlinks=true, citecolor = blue]{hyperref}

\usepackage{indentfirst}
\usepackage{changepage}
\usepackage{setspace}
\usepackage{indentfirst}
\usepackage{graphicx}
\usepackage{calc}

\allowdisplaybreaks

\theoremstyle{plain}
\newtheorem{theorem}{Theorem}[section]

\theoremstyle{remark}

\theoremstyle{definition}

\newenvironment{proof of theorem 1.1}{{\noindent \em Proof of Theorem 1.1.}}{\hfill $\Box$\par}
\newenvironment{proof of theorem 1.2}{{\noindent \em Proof of Theorem 1.2.}}{\hfill $\Box$\par}

\DeclareSymbolFont{EulerExtension}{U}{euex}{m}{n}
\DeclareMathSymbol{\euintop}{\mathop} {EulerExtension}{"52}
\DeclareMathSymbol{\euointop}{\mathop} {EulerExtension}{"48}

\begin{document}
	\title{Upper bounds for Erd\'{e}lyi's multivariate Laguerre polynomials}
	\author{Min-Jie Luo$^{\rm 1}$\thanks{Corresponding author}, Ravinder Krishna Raina$^{\rm 2}$}
	\date{}
	\maketitle
	\begin{center}\small
		$^{1}$\emph{Department of Mathematics, School of Mathematics and Statistics, \\ Donghua University, Shanghai 201620, \\ 
			People's Republic of China.}\\
		E-mail: \texttt{mathwinnie@live.com}, \texttt{mathwinnie@dhu.edu.cn}
	\end{center}
	\begin{center}\small
		$^{2}${\emph{M.P. University of Agriculture and Technology, Udaipur (Rajasthan), India\\
				\emph{Present address:} 10/11, Ganpati Vihar, Opposite Sector 5,\\
				Udaipur-313002, Rajasthan, India.}\\
			E-mail: \texttt{rkraina1944@gmail.com}; \texttt{rkraina\_7@hotmail.com}}  
	\end{center}
	
	
	\begin{abstract}

We establish in this paper two inequalities for the multivariate Laguerre polynomials introduced and studied by Arthur Erd\'{e}lyi [Sitzungsber. Akad. Wiss. Wien, Math.-Naturw. Kl., Abt. IIa 146 (1937), 431--467]. These inequalities generalize the well-known  Szeg\"{o}'s inequality for the Laguerre polynomials $L_n^{(\alpha)}(x)$. We also mention briefly few insightful remarks giving a comparative analysis concerning the upper bounds of the derived inequalities in the concluding section.\\

\\

\end{abstract}

\section{Introduction}\label{Introduction}

The Laguerre polynomials $L_n^{(\alpha)}(x)$ are usually defined by the generating function \cite[p. 449, Eq. (18.12.13)]{NIST Handbook}
\[
	(1-z)^{-\alpha-1}\exp\left(-\frac{xz}{1-z}\right)
	=\sum_{n=0}^{\infty}L_n^{(\alpha)}(x)z^n, ~ |z|<1.
\]
Explicitly, we have
\[
L_n^{(\alpha)}(x)=\frac{(1+\alpha)_n}{n!} {}_{1}F_{1}\left[\begin{matrix}
	-n\\
	\alpha+1
\end{matrix};x\right],
\]
where ${}_{1}F_{1}$ denotes the confluent hypergeometric function (see \cite[p. 443, Eq. (18.5.12)]{NIST Handbook}).

For the Laguerre polynomial $L_n^{(\alpha)}(x)$, the following inequality is well-known: 
\begin{equation}\label{SzegoIneq}
	|L_n^{(\alpha)}(x)|\leq\frac{(\alpha+1)_n}{n!}\mathrm{e}^{x/2},
\end{equation}
where $\alpha\geq 0$, $x\geq 0$ and $n\in\mathbb{Z}_{\geq0}:=\{0,1,2,\cdots\}$. Inequality \eqref{SzegoIneq} is usually called \emph{Szeg\"{o}'s inequality}. There are several different proofs of this famous result (see \cite{Ferreira-2025}, \cite{Petersen-Skovgaard-1952}, \cite{Szego-1918} and \cite{Watson-1939}). Several improvements have been suggested regarding the inequality \eqref{SzegoIneq}. Rooney \cite{Rooney-1984} extended the range of $\alpha$ to negative numbers, namely, 
\begin{equation}\label{RooneyInequality-1}
	|L_n^{(\alpha)}(x)|\leq2^{-\alpha}\mathrm{e}^{x/2},
\end{equation}
where $\alpha\leq 0$, $x\geq 0$ and  $n\in\mathbb{Z}_{\geq0}$. Later, Rooney in \cite{Rooney-1985} obtained the inequality that
\begin{equation}\label{RooneyInequality-2}
	|L_n^{(\alpha)}(x)|\leq q_n 2^{-\alpha}\mathrm{e}^{x/2},
\end{equation}
where $\alpha\leq -1/2$, $x\geq0$, $n\in\mathbb{Z}_{\geq0}$ and 
\begin{equation}\label{RooneyInequality-3}
	q_n:=\frac{((2n)!)^{1/2}}{2^{n+1/2}n!}\sim \frac{1}{\sqrt[4]{4\pi n}}~(n\rightarrow+\infty).
\end{equation}
As Rooney \cite{Rooney-1985}  has pointed out that the inequality \eqref{RooneyInequality-2} is stronger than \eqref{RooneyInequality-1} if $\alpha\leq -1/2$. For other inequalities for $L_n^{(\alpha)}(x)$ and the related functions, the interested readers may refer to \cite{Hardy-1939}, \cite{Lewandowski-Szynal-1998}, \cite{Love-1997}, \cite{Michalska-Szynal-2001} and \cite{Pogany-Srivastava-2007}.

There are many multivariate generalizations of the classical Laguerre polynomials in the literature (see, for example, \cite[p. 532, Eq. (1.4)]{Aktas-Erkus-Duman-2013}, \cite[p. 14, Eq. (1.1)]{Chak-1970}, \cite[p. 324, Definition 1]{Demenin-1971} and \cite[p. 387, Eq. (30)]{Liu-Lin-Lu-Srivastava-2013}). But in the present investigation, we focus ourselves on Erd\'{e}lyi's multivariate Laguerre polynomials \cite{Erdelyi-1937} which defines and introduces the Laguerre polynomials in $k$ variables by means of the multivariate generating function
\[
	(1-z_1-\cdots-z_k)^{-\alpha-1}
	\exp\left(-\frac{x_1 z_1+\cdots+x_k z_k}{1-z_1-\cdots-z_k}\right)
	=\sum_{n_1,\cdots,n_k=0}^{\infty}L_{n_1,\cdots,n_k}^{(\alpha)}(x_1,\cdots,x_k)z_1^{n_1}\cdots z_k^{n_k},
\]
where $|z_1|+\cdots+|z_k|<1$. For convenience, we let $\mathbf{x}=(x_1,\cdots,x_k)$ so that $\lambda \mathbf{x}=(\lambda x_1,\cdots,\lambda x_k)$ $(\lambda\in\mathbb{R})$. Erd\'{e}lyi \cite[p. 458, Eq. (11.3)]{Erdelyi-1937} has shown that 
\[
L_{n_1,\cdots,n_k}^{(\alpha)}(\mathbf{x})
=\frac{(\alpha+1)_{n_1+\cdots+n_k}}{n_1!\cdots n_k!}\Phi_2^{(k)}\left[-n_1,\cdots,-n_k;\alpha+1;\mathbf{x}\right],
\]
where $\Phi_2^{(k)}$ is defined by (see \cite[p. 446]{Erdelyi-1937} and \cite[p. 34]{Srivastava-Karlsson-Book-1985})
\begin{equation}\label{ConfluentLauricella}
	\Phi_2^{(k)}\left[b_1,\cdots,b_k;c;\mathbf{x}\right]
	:=\sum_{j_1=0}^{\infty}\cdots\sum_{j_k=0}^{\infty}
	\frac{(b_1)_{j_1}\cdots (b_k)_{j_k}}{(c)_{j_1+\cdots+j_k}}\frac{x_1^{j_1}}{j_1!}\cdots \frac{x_k^{j_k}}{j_k!}.
\end{equation}

It may be pointed out here that Erd\'{e}lyi also obtained
an elegant generalization of \emph{Hardy-Hille formula} in \cite {Erdelyi-1937} for $L_{n_1,\cdots,n_k}^{(\alpha)}(\mathbf{x})$ which in our opinion makes this set of polynomials particularly useful and of importance. Carlitz \cite{Carlitz-1970} further provided an elementary proof of this formula and considered other new generalizations.

Our main result is contained in the following theorem.
\begin{theorem}\label{MainTheorem}
	Let $\alpha>0$, $\mathbf{x}=(x_1,\cdots,x_k)$ $(x_j\geq0,~j=1,\cdots,k)$ and $\displaystyle \|\mathbf{x}\|:=\max_{1\leq j\leq k}\{|x_j|\}$. Then 
	\begin{equation}\label{MainTheorem-1}
		\left|L_{n_1,\cdots,n_k}^{(\alpha)}(\mathbf{x})\right|
		\leq 
		2^{k-1} \frac{(\alpha+1)_{n_1+\cdots+n_k}}{(\frac{1}{k})_{n_1}\cdots (\frac{1}{k})_{n_k}}\mathrm{e}^{\|\mathbf{x}\|/2}.
	\end{equation}
\end{theorem}

When $k=1$, \eqref{MainTheorem-1} reduces to Szeg\"{o}'s inequality \eqref{SzegoIneq}. Another variation of \eqref{MainTheorem-1} can be considered which holds for certain extended range of parameter $\alpha$. This result is given by the following theorem.

\begin{theorem}\label{MainTheorem-B}
	Let $\alpha>-1$, $\mathbf{x}=(x_1,\cdots,x_k)$ $(x_j\geq0,~j=1,\cdots,k)$ and $\displaystyle \|\mathbf{x}\|:=\max_{1\leq j\leq k}\{|x_j|\}$. Then 
	\begin{equation}\label{MainTheorem-B-1}
		\left|L_{n_1,\cdots,n_k}^{(\alpha)}(\mathbf{x})\right|
		\leq q_{n_1}\cdots q_{n_k} 2^{k-\frac{1}{2}} 
		\frac{(\alpha+1)_{n_1+\cdots+n_k}}{(\frac{1}{2k})_{n_1}\cdots (\frac{1}{2k})_{n_k}}
		\mathrm{e}^{\|\mathbf{x}\|/2},
	\end{equation}
	where $q_n$ is given by \eqref{RooneyInequality-3}. 
\end{theorem}

The proofs of Theorem \ref{MainTheorem} and Theorem \ref{MainTheorem-B} will be presented in Section \ref{Proof of Th 1.1} and Section \ref{Proof of Th 1.2}, respectively.
A comparative analysis of the upper bounds for the multivariate Laguerre polynomials defined above along with some insightful remarks are mentioned in Section \ref{Concluding remarks}.

\section{Proof of Theorem \ref{MainTheorem}}\label{Proof of Th 1.1}

Our starting point is an integral representation due to Srivastava and Niukkanen \cite[p. 249, Eq. (22)]{Srivastava-Niukkanen}. To state their result in a compact form, we introduce the concept of Dirichlet measure \cite[p. 64, Definition 4.4-1]{Carlson-Book-1977}. Let $E_k$ be the standard simplex in $\mathbb{R}^k$. The Dirichlet measure $\mu_{(b_1,\cdots,b_k,\beta)}(\mathbf{u})$ is defined on $E$ by
\begin{equation}\label{DirichletMeasure}
	\mathrm{d}\mu_{(b_1,\cdots,b_k,\beta)}(\mathbf{u})
	=\frac{\Gamma(b_1+\cdots+b_{k}+\beta)}{\Gamma(b_1)\cdots\Gamma(\beta_{k})\Gamma(\beta)}
	u_1^{b_1-1}\cdots u_k^{b_k-1}(1-u_1-\cdots-u_k)^{\beta-1}\mathrm{d}u_1\cdots\mathrm{d}u_k,
\end{equation}
where $b_j>0$ $(j=1,\cdots,k)$ and $\beta>0$. We have
\[
\int\cdots\int_{E_k} \mathrm{d}\mu_{(b_1,\cdots,b_k,\beta)}(\mathbf{u})=1.
\]
The integral representation given by Srivastava and Niukkanen \cite{Srivastava-Niukkanen} can now be used for the multivariable Laguerre polynomials, and in view of \eqref{DirichletMeasure}, we can express
\begin{align}\label{Srivastava-Niukkanen Integral}
	L_{n_1,\cdots,n_k}^{(\alpha_1+\cdots+\alpha_k+\beta+k)}(\mathbf{x})
	&=\frac{(\alpha_1+\cdots+\alpha_k+\beta+k+1)_{n_1+\cdots+n_k}}{(\alpha_1+1)_{n_1}\cdots (\alpha_k+1)_{n_k}}\notag\\
	&\cdot 
	\int\cdots\int_{E_k} 
	\prod_{j=1}^{k} 
	L_{n_j}^{(\alpha_j)}(x_j u_j)\mathrm{d}\mu_{(\alpha_1+1,\cdots,\alpha_k+1,\beta+1)}(\mathbf{u}),
\end{align}
where $\alpha_j>-1$ $(j=1,\cdots,k)$ and $\beta>-1$.

If we set $\alpha_1=\cdots=\alpha_k=(1-k)/k\in(-1,0]$ and $\beta=\alpha-1$ $(\alpha>0)$ in \eqref{Srivastava-Niukkanen Integral}, then we get 
\[
	L_{n_1,\cdots,n_k}^{(\alpha)}(\mathbf{x})
	=\frac{(\alpha+1)_{n_1+\cdots+n_k}}{(\frac{1}{k})_{n_1}\cdots (\frac{1}{k})_{n_k}}
	\int\cdots\int_{E_k} 
	\prod_{j=1}^{k} 
	L_{n_j}^{(\frac{1-k}{k})}(x_j u_j)\mathrm{d}\mu_{(\frac{1}{k},\cdots,\frac{1}{k},\alpha)}(\mathbf{u}).
\]
Making use of Rooney's inequality \eqref{RooneyInequality-1}, we have
\begin{align}\label{MainTheorem-Proof-1}
	\left|L_{n_1,\cdots,n_k}^{(\alpha)}(\mathbf{x})\right| 
	&\leq \frac{(\alpha+1)_{n_1+\cdots+n_k}}{(\frac{1}{k})_{n_1}\cdots (\frac{1}{k})_{n_k}}
	\int\cdots\int_{E_k} 
	\prod_{j=1}^{k} 
	\left|L_{n_j}^{(\frac{1-k}{k})}(x_j u_j)\right|
	\mathrm{d}\mu_{(\frac{1}{k},\cdots,\frac{1}{k},\alpha)}(\mathbf{u})\notag\\
	&\leq \frac{(\alpha+1)_{n_1+\cdots+n_k}}{(\frac{1}{k})_{n_1}\cdots (\frac{1}{k})_{n_k}}\cdot 2^{k-1}
	\int\cdots\int_{E_k}
	\mathrm{e}^{(x_1 u_1+\cdots+x_k u_k)/2}
	\mathrm{d}\mu_{(\frac{1}{k},\cdots,\frac{1}{k},\alpha)}(\mathbf{u}),
\end{align}            
where $\alpha>0$.

Note that
\begin{align}
	&\int\cdots\int_{E_k}
	\mathrm{e}^{(x_1 u_1+\cdots+x_k u_k)/2}
	\mathrm{d}\mu_{(\frac{1}{k},\cdots,\frac{1}{k},\alpha)}(\mathbf{u})\notag \\
	&\hspace{0.5cm}=\sum_{j_1=0}^{\infty}\cdots\sum_{j_k=0}^{\infty}\frac{(\frac{1}{2}x_1)^{j_1}}{j_1!}\cdots\frac{(\frac{1}{2}x_k)^{j_k}}{j_k!}\int\cdots\int_{E_k}
	u_1^{j_1}\cdots u_k^{j_k}
	\mathrm{d}\mu_{(\frac{1}{k},\cdots,\frac{1}{k},\alpha)}(\mathbf{u})\notag \\ 
	&\hspace{0.5cm}=\sum_{j_1=0}^{\infty}\cdots\sum_{j_k=0}^{\infty}
	\frac{(\frac{1}{k})_{j_1}\cdots(\frac{1}{k})_{j_k}}{(\alpha+1)_{j_1+\cdots+j_k}}\frac{(\frac{1}{2}x_1)^{j_1}}{j_1!}\cdots\frac{(\frac{1}{2}x_k)^{j_k}}{j_k!}\notag \\ 
	&\hspace{0.5cm}=\Phi_2^{(k)}\left[\frac{1}{k},\cdots,\frac{1}{k};\alpha+1;\frac{1}{2}\mathbf{x}\right], \label{MainTheorem-Proof-2}
\end{align}
where $\Phi_2^{(k)}$ is defined by \eqref{ConfluentLauricella}. Thus, it follows from \eqref{MainTheorem-Proof-1} and \eqref{MainTheorem-Proof-2} that
\begin{equation}\label{MainTheorem-Proof-3}
	\left|L_{n_1,\cdots,n_k}^{(\alpha)}(\mathbf{x})\right|
	\leq
	2^{k-1} \frac{(\alpha+1)_{n_1+\cdots+n_k}}{(\frac{1}{k})_{n_1}\cdots (\frac{1}{k})_{n_k}} \Phi_2^{(k)}\left[\frac{1}{k},\cdots,\frac{1}{k};\alpha+1;\frac{1}{2}\mathbf{x}\right].
\end{equation}

In order to simplify the upper bound in \eqref{MainTheorem-Proof-3}, we note that
\[
\Phi_2^{(k)}\left[\frac{1}{k},\cdots,\frac{1}{k};\alpha+1;\frac{1}{2}\mathbf{x}\right]
\leq \sum_{j_1=0}^{\infty}\cdots\sum_{j_k=0}^{\infty}
\frac{(\frac{1}{k})_{j_1}\cdots(\frac{1}{k})_{j_k}}{(\alpha+1)_{j_1+\cdots+j_k}}\frac{(\frac{1}{2}\|\mathbf{x}\|)^{j_1+\cdots+j_k}}{j_1!\cdots j_k!}.
\]
Now to reduce the multiple series into a single series, we use the following identity due to Panda \cite[p. 166, Theorem 2]{Panda-1974}:
\[
	\sum_{j_1=0}^{\infty}\cdots \sum_{j_k=0}^{\infty}
	C(j_1+\cdots+j_k)(\alpha_1)_{j_1}\cdots(\alpha_k)_{j_k}\frac{x^{j_1+\cdots+j_k}}{j_1!\cdots j_k!}=\sum_{j=0}^{\infty}C(j)(\alpha_1+\cdots+\alpha_k)_j\frac{x^j}{j!},
\]
and therefore
\begin{equation}\label{MainTheorem-Proof-5}
	\Phi_2^{(k)}\left[\frac{1}{k},\cdots,\frac{1}{k};\alpha+1;\frac{1}{2}\mathbf{x}\right]
	\leq \sum_{j=0}^{\infty}
	\frac{(1)_{j}}{(\alpha+1)_{j}}\frac{(\frac{1}{2}\|\mathbf{x}\|)^{j}}{j!}
	={}_{1}F_{1}\left[\begin{matrix}
		1\\
		\alpha+1
	\end{matrix};\frac{1}{2}\|\mathbf{x}\|\right]
	\leq \mathrm{e}^{\|\mathbf{x}\|/2}.
\end{equation}

Combining now \eqref{MainTheorem-Proof-3} and \eqref{MainTheorem-Proof-5},  we obtain the desired inequality \eqref{MainTheorem-1}.
This completes the proof.

\section{Proof of Theorem \ref{MainTheorem-B}}\label{Proof of Th 1.2}

Let us set
$\alpha_1=\cdots=\alpha_k=(1-2k)/(2k)\in(-1,-1/2]$
and $\beta=\alpha-(1/2)$ $(\alpha>-1/2)$ in \eqref{Srivastava-Niukkanen Integral} so that
\[
L_{n_1,\cdots,n_k}^{(\alpha)}(\mathbf{x})
=\frac{(\alpha+1)_{n_1+\cdots+n_k}}{(\frac{1}{2k})_{n_1}\cdots (\frac{1}{2k})_{n_k}}
\int\cdots\int_{E_k} 
\prod_{j=1}^{k} 
L_{n_j}^{(\frac{1-2k}{2k})}(x_j u_j)\mathrm{d}\mu_{(\frac{1}{2k},\cdots,\frac{1}{2k},\alpha+\frac{1}{2})}(\mathbf{u}).
\]
Then, making use of Rooney's inequality \eqref{RooneyInequality-2}, we have 
\begin{align*}
	\left|L_{n_1,\cdots,n_k}^{(\alpha)}(\mathbf{x})\right| 
	&\leq \frac{(\alpha+1)_{n_1+\cdots+n_k}}{(\frac{1}{2k})_{n_1}\cdots (\frac{1}{2k})_{n_k}}
	\int\cdots\int_{E_k} 
	\prod_{j=1}^{k} 
	\left|L_{n_j}^{(\frac{1-2k}{2k})}(x_j u_j)\right|
	\mathrm{d}\mu_{(\frac{1}{2k},\cdots,\frac{1}{2k},\alpha+\frac{1}{2})}(\mathbf{u})\\
	&\leq \frac{(\alpha+1)_{n_1+\cdots+n_k}}{(\frac{1}{2k})_{n_1}\cdots (\frac{1}{2k})_{n_k}}\cdot q_{n_1}\cdots q_{n_k} 2^{\frac{2k-1}{2}}
	\int\cdots\int_{E_k}
	\mathrm{e}^{(x_1 u_1+\cdots+x_k u_k)/2}
	\mathrm{d}\mu_{(\frac{1}{2k},\cdots,\frac{1}{2k},\alpha+\frac{1}{2})}(\mathbf{u}),
\end{align*}
where $\alpha>-1/2$.

Note that
\begin{align*}
	&\int\cdots\int_{E_k}
	\mathrm{e}^{(x_1 u_1+\cdots+x_k u_k)/2}
	\mathrm{d}\mu_{(\frac{1}{2k},\cdots,\frac{1}{2k},\alpha+\frac{1}{2})}(\mathbf{u})\\
	&\hspace{0.5cm}=\sum_{j_1=0}^{\infty}\cdots\sum_{j_k=0}^{\infty}\frac{(\frac{1}{2}x_1)^{j_1}}{j_1!}\cdots\frac{(\frac{1}{2}x_k)^{j_k}}{j_k!}\int\cdots\int_{E_k}
	u_1^{j_1}\cdots u_k^{j_k}
	\mathrm{d}\mu_{(\frac{1}{2k},\cdots,\frac{1}{2k},\alpha+\frac{1}{2})}(\mathbf{u})\\
	&\hspace{0.5cm}=\sum_{j_1=0}^{\infty}\cdots\sum_{j_k=0}^{\infty}
	\frac{(\frac{1}{2k})_{j_1}\cdots(\frac{1}{2k})_{j_k}}{(\alpha+1)_{j_1+\cdots+j_k}}\frac{(\frac{1}{2}x_1)^{j_1}}{j_1!}\cdots\frac{(\frac{1}{2}x_k)^{j_k}}{j_k!}\\
	&\hspace{0.5cm}=\Phi_2^{(k)}\left[\frac{1}{2k},\cdots,\frac{1}{2k};\alpha+1;\frac{1}{2}\mathbf{x}\right],
\end{align*}
where $\Phi_2^{(k)}$ is defined by \eqref{ConfluentLauricella}. Thus
\begin{equation}\label{MainTheorem-B-Proof-2}
	\left|L_{n_1,\cdots,n_k}^{(\alpha)}(\mathbf{x})\right|
	\leq \frac{(\alpha+1)_{n_1+\cdots+n_k}}{(\frac{1}{2k})_{n_1}\cdots (\frac{1}{2k})_{n_k}}\cdot q_{n_1}\cdots q_{n_k} 2^{\frac{2k-1}{2}}\Phi_2^{(k)}\left[\frac{1}{2k},\cdots,\frac{1}{2k};\alpha+1;\frac{1}{2}\mathbf{x}\right].
\end{equation}

Following similar approach as in the proof of Theorem \ref{MainTheorem} in Section \ref{Proof of Th 1.1} above, the function $\Phi_2^{(k)}\left[\frac{1}{2k},\cdots,\frac{1}{2k};\alpha+1;\frac{1}{2}\mathbf{x}\right] $ occurring in the upper bound of \eqref{MainTheorem-B-Proof-2} simplifies to
\begin{align}\label{MainTheorem-B-Proof-3}
	\Phi_2^{(k)}\left[\frac{1}{2k},\cdots,\frac{1}{2k};\alpha+1;\frac{1}{2}\mathbf{x}\right]
	&\leq \sum_{j_1=0}^{\infty}\cdots\sum_{j_k=0}^{\infty}
	\frac{(\frac{1}{2k})_{j_1}\cdots(\frac{1}{2k})_{j_k}}{(\alpha+1)_{j_1+\cdots+j_k}}\frac{(\frac{1}{2}\|\mathbf{x}\|)^{j_1+\cdots+j_k}}{j_1!\cdots j_k!}\notag \\
	&=\sum_{j=0}^{\infty}
	\frac{(\frac{1}{2})_{j}}{(\alpha+1)_{j}}\frac{(\frac{1}{2}\|\mathbf{x}\|)^{j}}{j!}
	={}_{1}F_{1}\left[\begin{matrix}
		\frac{1}{2}\\
		\alpha+1
	\end{matrix};\frac{1}{2}\|\mathbf{x}\|\right]
	\leq \mathrm{e}^{\|\mathbf{x}\|/2}.                                               
\end{align}
Using \eqref{MainTheorem-B-Proof-2} and \eqref{MainTheorem-B-Proof-3}, we obtain the desired inequality \eqref{MainTheorem-B-1}. This completes the proof.

\section{Concluding remarks}\label{Concluding remarks}

\begin{itemize}
	\item[(i)] We have mentioned in Section \ref{Introduction} that Szeg\"{o}'s inequaltiy \eqref{SzegoIneq} is reproduced when taking $k=1$ in \eqref{MainTheorem-1}. However, letting $k=1$ in inequality \eqref{MainTheorem-B-1} of Theorem \ref{MainTheorem-B} gives
	\[
	\left|L_{n}^{(\alpha)}(x)\right|
	\leq \Bigg(\frac{(1)_n}{(\frac{1}{2})_n}\Bigg)^{1/2}
	\frac{(\alpha+1)_{n}}{n!}
	\mathrm{e}^{x/2}.
	\]
	The additional factor $((1)_n/(\tfrac{1}{2})_n)^{1/2}$ makes the upper bound a little worse than \eqref{SzegoIneq} when $\alpha>0$. But we can show that, for $k\geq2$, inequality \eqref{MainTheorem-B-1} could be better than \eqref{MainTheorem-1}. For this purpose, we let $n_1=\cdots=n_k=n$ in Theorem \ref{MainTheorem} and Theorem \ref{MainTheorem-B}. The resulting inequalities are
	\[
	\left|L_{n,\cdots,n}^{(\alpha)}(\mathbf{x})\right|
	\leq 
	\mathsf{A}_n(\alpha,k)\mathrm{e}^{\|\mathbf{x}\|/2}
	~~~
	\text{and} 
	~~~
	\left|L_{n,\cdots,n}^{(\alpha)}(\mathbf{x})\right|
	\leq \mathsf{B}_n(\alpha,k)
	\mathrm{e}^{\|\mathbf{x}\|/2},
	\] 
	where
	\[
	\mathsf{A}_n(\alpha,k):= 2^{k-1} \frac{(\alpha+1)_{k n}}{((\frac{1}{k})_{n})^k}~~~\text{and}~~~
	\mathsf{B}_n(\alpha,k):=(q_{n})^k 2^{k-\frac{1}{2}} 
	\frac{(\alpha+1)_{kn}}{((\frac{1}{2k})_{n})^k}.
	\]
	We have after a little computation that
	\[
	\frac{\mathsf{A}_n(\alpha,k)}{\mathsf{B}_n(\alpha,k)}\sim
	2^{\frac{k-1}{2}}\pi^{\frac{k}{4}} \left(\frac{\Gamma(\frac{1}{2k})}{\Gamma(\frac{1}{k})}\right)^k n^{\frac{k}{4}-\frac{1}{2}}~(n\rightarrow+\infty).
	\]
	It is therefore clear that when $k\geq 2$, the upper bound provided by \eqref{MainTheorem-B-1} would be better than the one given by \eqref{MainTheorem-1}. In addition, if we want to find an accurate estimate of $L_{n,\cdots,n}^{(\alpha)}(\mathbf{x})$, it would be possible to apply the so-called \emph{diagonal method} (see \cite[Section 13.1]{Pemantle-Wilson-Book-2013}).  We may have to first find a closed form for the \emph{diagonal generating function}
	\begin{equation}\label{DiagonalGF}
	\sum_{n=0}^{\infty}L_{n,\cdots,n}^{(\alpha)}(\mathbf{x}) z^n.
	\end{equation}
	It seems at this point very difficult to find such a closed-form result for \eqref{DiagonalGF} for $k\geq 3$, and therefore we leave it for future research. 
	
	\item[(ii)] Our main technique employed in the present investigation was to make the parameters 
	$\alpha_1, \cdots,\alpha_k$ in the integral representation \eqref{Srivastava-Niukkanen Integral} sufficiently close to $-1$. Based on this observation, it seems necessary and worthwhile to mention here the following inequality discovered by Lewandowski and Szynal \cite[p. 532]{Lewandowski-Szynal-1998}: 
	\begin{equation}\label{Lewandowski-Szynal Inequality}
		|L_n^{(\alpha)}(x)|\leq\frac{(\alpha+1)_n}{n!}\sigma_n^{(\alpha)}(\mathrm{e}^x),
	\end{equation}
	where $\alpha\geq-1/2$, $x\geq0$,  $n\in\mathbb{Z}_{\geq0}$ and
	\[
	\sigma_{n}^{(\alpha)}(\mathrm{e}^x)=\frac{n!}{(\alpha+1)_n}\sum_{k=0}^{n}\frac{(\alpha+1)_{n-k}}{(n-k)!}\frac{x^k}{k!}.
	\]
	This inequality \eqref{Lewandowski-Szynal Inequality} is better than Szeg\"{o}'s inequality \eqref{SzegoIneq} for large $x$, because $\sigma_{n}^{(\alpha)}(\mathrm{e}^x)$ is a polynomial. Inequality \eqref{Lewandowski-Szynal Inequality} has a remarkable application. By using \eqref{Lewandowski-Szynal Inequality}, Luo and Raina \cite{Luo-Raina-2020} established an inequality for the associated Pollaczek polynomials.  Since the inequality \eqref{Lewandowski-Szynal Inequality} holds for all $\alpha\geq -1/2$, it can therefore also be used to obtain other inequalities for the multivariate Laguerre polynomials. We choose to omit such results here because derivations of these results would not involve any new ideas in their proofs. 
\end{itemize}

\section*{Acknowledgement}

The research of the first author is supported by National Natural Science Foundation of China (Grant No. 12001095).

\end{document}